\documentclass[preprint,review,11pt]{elsarticle}

\newcommand\BibTeX{{\rmfamily B\kern-.05em \textsc{i\kern-.025em b}\kern-.08em
T\kern-.1667em\lower.7ex\hbox{E}\kern-.125emX}}

\usepackage{etex} 
\usepackage{graphicx,subfigure,color,colordvi}
\usepackage{latexsym,amsbsy,epsfig,fancybox}
\usepackage{amstext}
\usepackage{amsfonts,textgreek}
\usepackage{mathrsfs}
\usepackage{mathtools}
\usepackage{stmaryrd,dsfont}
\usepackage{bbm}
\usepackage{url}
\usepackage{wrapfig,bm}
\usepackage{tikz}
\usepackage{pstricks}
\usepackage{caption}
\usepackage{psfrag}
\usepackage{upgreek}
\usepackage{isomath}
\usepackage{latexsym}
\usepackage{graphicx}
\usepackage{array}
\usepackage{multirow}
\usepackage{booktabs}
\usepackage{colortbl}
\usepackage{verbatim}
\usepackage{epstopdf,overpic}
\usepackage{enumerate}
\usepackage{amsmath,amssymb,xspace}

\usepackage{pgfplots}
\usetikzlibrary{plotmarks}

\newcommand{\eref}[1]{Equation~(\ref{#1})}

\newcommand{\fref}[1]{Figure~\ref{#1}}

\newcommand{\vm}[1]{\bm{\mathrm{#1}}} 
\newcommand{\bsym}[1]{\bm{#1}}

\newcommand{\mat}[1]{\bm{\mathrm{#1}}} 
\newcommand{\transpose}{\mathrm{T}}

\tikzstyle{nicebox}=[draw=black!100, fill=white!10, rectangle, inner sep=4pt, inner ysep=16pt]
\tikzstyle{niceboxtitle}=[draw=black!100, fill=white, text=black, rectangle]

\usetikzlibrary{arrows,decorations.markings}

\definecolor{forestgreen}{RGB}{34, 139, 34}
\definecolor{lightgray}{gray}{0.92}

\newtheorem{remark}{Remark}

\newcommand{\xx}{\vm{x}}
\newcommand{\rmd}{\mathrm{d}}

\newcommand{\uu}{\mat{u}}

\journal{arXiv}

\begin{document}

\begin{frontmatter}

\title{Linear smoothed extended finite element method}

\author[India1,India3]{M Surendran}
\author[India2]{S Natarajan}
\author[lux]{SPA Bordas \corref{cor1}\fnref{fn1}}
 \author[India1,India3]{GS Palani}
 
\address[India1]{Academy of Scientific and Innovative Research, CSIR-Structural Engineering Research Center, CSIR Campus, Chennai - 600107, India.}
\address[India3]{Center for Scientific and Industrial Research, Chennai - 600107, India.}
\address[India2]{Department of Mechanical Engineering, Indian Institute of Technology, Madras, Chennai - 600036, India.}
\address[lux]{Facult\'e des Sciences, de la Technologie et de al Communication, University of Luxembourg, Luxembourg.}

\fntext[fn1]{Facult\'e des Sciences, de la Technologie et de al Communication, University of Luxembourg, Luxembourg., Email: stephane.bordas@uni.lu}

\begin{abstract}
The extended finite element method (XFEM) was introduced in 1999 to treat problems involving discontinuities with no or minimal remeshing through appropriate enrichment functions. This enables elements to be split by a discontinuity, strong or weak and hence requires the integration of discontinuous functions or functions with discontinuous derivatives over elementary volumes. Moreover, in the case of open surfaces and singularities, special, usually non-polynomial functions must also be integrated. A variety of approaches have been proposed to facilitate these special types of  numerical integration, which have been shown to have a large impact on the accuracy and the convergence of the numerical solution. The smoothed extended finite element method (SmXFEM) \cite{bordas2010strain}, for example, makes numerical integration elegant and simple by transforming volume integrals into surface integrals. However, it was reported in \cite{bordas2010strain, bordas2011performance} that the strain smoothing is inaccurate when non-polynomial functions are in the basis. This is due to the constant smoothing function used over the smoothing domains which destroys the effect of the singularity. In this paper, we investigate the benefits of a recently developed Linear smoothing procedure \cite{francis2015linear} which provides better approximation to higher order polynomial fields in the basis. Some benchmark problems in the context of linear elastic fracture mechanics (LEFM) are solved to compare the standard XFEM, the constant-smoothed XFEM (Sm-XFEM) and the linear-smoothed XFEM (LSm-XFEM). We observe that the convergence rates of all three methods are the same. The stress intensity factors (SIFs) computed through the proposed LSm-XFEM are however more accurate than that obtained through Sm-XFEM. To conclude, compared to the conventional XFEM, the same order of accuracy is achieved at a relatively low computational effort.
\end{abstract}


\end{frontmatter}

\vspace{-20pt}


\section{Introduction}
\label{intro}

The major difficulty associated with solving problems involving evolving discontinuities is the meshing and re-meshing required as the discontinuities evolve in time. This difficulty is exacerbated when singularities are also present, as is the case in crack growth simulations. These difficulties are somewhat alleviated by the introduction of enrichment functions to represent the discontinuities and the singularities at the patch level, in finite elements or meshfree methods.  A first approach to treat discontinuities without an explicit meshing was proposed as early as 1995 in \cite{oliyer1995continuum}. A much more versatile approach was presented a few years later in the form of the extended finite element method (XFEM) \cite{belytschko1999elastic} \cite{dolbow1999finite} by exploiting the partition of unity property identified by Melenk and Babu\v ska \cite{melenk1996partition}. Partition of unity enrichment for problems with discontinuous solutions is now widely used both in academia and in industrial practice \cite{bordasmoran2006,bordasconley2007,wyartcoulon2007} and is known under various names, including the generalized finite element method (GFEM) \cite{duartehamzeh2001,babuskabanerjee2012} and the extended finite element method (XFEM). The approach has also been widely used in the form of enriched meshfree methods \cite{rabczukbordas2007}. 

Another problem associated with partition of unity methods involving non-polynomial basis functions is to integrate the resulting fields accurately. These enriched methods, also carry along the element mapping involved in building the system matrices. The regularity and positive definiteness of the isoparametric mapping poses a number of restrictions on the allowable shapes of the finite elements: for example, the element should be convex. Meshfree methods also have to  face such problems associated with the regularity, distortion and clustering in the point cloud. Under large distortions, meshfree methods face numerical instabilities and low accuracy \cite{belytschkoguo2000}.  Nodal integration also leads to instabilities in cases where higher order shape functions are used. This is due to the fact that in the meshfree methods each node would be associated with a support domain. And the shape functions derivatives would be integrated in this support domain. This implies that each integration domain would be associated with only one integration point (i.e the node). Hence when only one integration is point is considered for higher order functions (other than constant strain) the integration would be similar to the inadequate reduced integration which in turn causes instabilities.

To alleviate these instabilities, the strain smoothing concept was introduced for  meshfree methods \cite{chen2001stabilized}. The basic idea of strain smoothing is to transform numerical integration over volumes to integration over surfaces, thereby removing instabilities due to the evaluation of the shape function derivatives at the nodes. This approach was later extended to finite element methods by Liu et al \cite{liu2007smoothed}. The resulting method was coined the Smoothed finite element method (SFEM), was discussed in a number of papers \cite{bordas2010approximation,bordas2011performance,liu2009edge,liu2009node,nguyenrabczuk2008,nguyen2008smooth} and applied to a wide variety of problems. The SFEM reduces the mesh sensitivity to some extent by avoiding the necessity of evaluating the Jacobian. Since the derivatives are not needed, the iso-parametric mapping is also avoided.

The SFEM involves computation of a smoothed strain from the standard compatible strain field. The smoothed strain is evaluated as a spatial average of the standard strain field over smoothing subcells which cover the domain andthat can be chosen independently from the mesh structure. These smoothing cells are typically constructed from the mesh following different approaches. This gives rise to a number of methods including cell-based SFEM (CS-FEM) \cite{liu2007smoothed, dai2007n, nguyen2008smooth, bordas2010approximation}, node-based SFEM (NS-FEM) \cite{liu2009node}, edge-based SFEM (ESFEM) \cite{liu2009edge}, face-based SFEM (FS-FEM). 

The method was later extended to solve problems with field discontinuities, both strong and weak, by Bordas et al \cite{bordas2010strain}. This was achieved by extending strain smoothing to the partition of unity framework \cite{babuvska1997partition, melenk1996partition}. Though the smoothed FEM did make the integration procedure more elegant, it was also reported in \cite{bordas2011performance} that the error levels are higher due to the inaccurate approximation of the near tip singular fields. Similar errors were also encountered in higher order elements and polygons \cite{natarajan2014convergence}. It is noteworthy that similar difficulties are also present in meshfree methods, as  addressed in \cite{duan2012second} by employing the Discrete Divergence Consistency (DDC) requirement by establishing a compatibility relation between the shape function and its derivatives. This was recently extended for the FEM in \cite{francis2015linear} and named: Linear smoothing (LS) scheme. It was also reported that the linear smoothing scheme provides an improved accuracy compared to the standard constant-based smoothing of the SFEM. 

The present paper aims at investigating how the use of higher-order smoothing, in particular linear smoothing, resolves the lack of accuracy observed when constant smoothing is employed with non-polynomial bases functions. The paper is organized as follows. After presenting the governing equations for elasto-statics, a brief overview of the extended finite element method is given in Section \ref{xfem}. Section \ref{ssm} presents the linear smoothing technique. A few standard benchmark problems in linear elastic fracture mechanics, solved by using the developed method are presented and the accuracy, convergence and the efficiency of the proposed method are discussed in Section \ref{numerics}, followed by concluding remarks in the last section.



\section{Theoretical formulation}
\label{xfem}
\subsection{Governing equations for elastostatics}
\begin{figure}[h]
\centering
\includegraphics[width=0.7\textwidth]{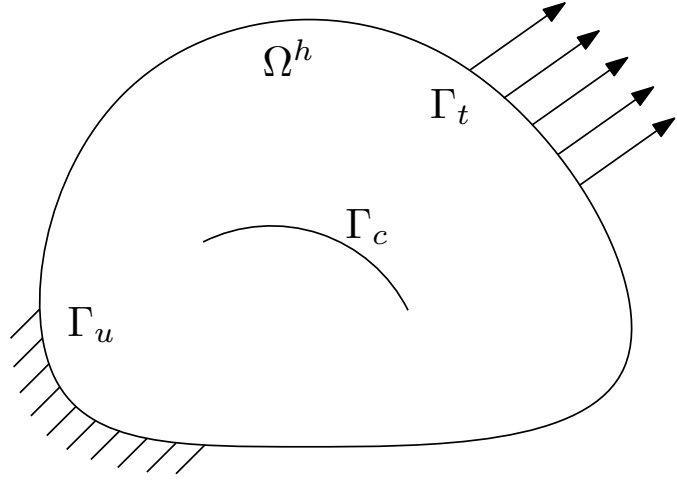}
\caption{General elastic body with an internal discontinuity, Neumann and Dirichlet boundary conditions.}
\label{body}
\end{figure}

Consider a linear elastic body as shown in \fref{body}, with a discontinuity. Let the domain be $\Omega \in \mathbb R^d $ bounded by $\Gamma$. In this case the boundary has three parts namely $\Gamma_u$, where the displacement boundary conditions are applied, $\Gamma_t$, where the tractive boundary conditions are applied and $\Gamma_c$, which is the traction free surface representing the discontinuity, such that $\Gamma=\Gamma_u \cup \Gamma_t \cup \Gamma_c$ and $ \Gamma_u \cap \Gamma_t = \emptyset $. 
The governing equation to be solved is
\begin{equation}
\bm{\nabla} \cdot \bm{\sigma} +\mathbf{b}= \mathbf{0} \quad\textrm{in } \Omega
\label{eq14}
\end{equation}
The boundary conditions are as follows
\begin{subequations}
\begin{equation}
\bm{\sigma} \cdot \mathbf{n} = \mathbf{0} \quad\textrm{on } \Gamma_c
\label{eq15}
\end{equation}
\begin{equation}
\bm{\sigma} \cdot \mathbf{n}=\mathbf{\hat{t}} \quad\textrm{on } \Gamma_t
\label{eq16a}
\end{equation}
\begin{equation}
\mathbf{u=\bar{u}} \quad\textrm{on } \Gamma_u
\label{eq16b}
\end{equation}
\label{eq16}
\end{subequations}
where, is the gradient operator, $\bm{\sigma}$ is the Cauchy stress tensor, $\mathbf{b}$ is the body force per unit volume, $\bf{n}$ is the unit outward normal and $\mathbf{t}$ is the applied tractive stress. For a body undergoing small displacements and strains, the strain displacement equation reads as
\begin{equation}
\bm{\varepsilon}=\bm{\nabla}_s \bf{u}
\label{eq17}
\end{equation}
where, $\nabla_s$ is the symmetric part of the gradient operator. By substituting the constitutive relations and the strain-displacement relations the weak form of the above Equation \eqref{eq14} becomes Equation\eqref{eq18} in the absence of the body forces: find $\uu \in \mathscr{U}$ such that
\begin{equation}
\int\limits_\Omega \boldsymbol{\varepsilon}(\mathbf{u}) \colon \mathbf{C} \colon \boldsymbol{\varepsilon}(\mathbf{v})~\mathrm{d}\Omega = \int\limits_{\Gamma_t} \mathbf{t} \mathbf{v}~\mathrm{d}\Gamma
\label{eq18}
\end{equation}
where, $\bf{u}$ and $\bf{v}$ are the trial and the test function spaces, respectively. For the aforementioned problem, the function spaces includes functions that are discontinuous across $\Gamma_c$.
\begin{subequations}
\begin{align}
\mathscr{U} &:=
\left\{\vm{u}(\vm{x})\in [ C^0(\Omega)]^d : \vm{u} \in
[ \mathcal{W}(\Omega)]^d \subseteq [ H^{1}(\Omega)]^d, \ \vm{u} = \hat{\vm{u}}
\ \textrm{on } \Gamma_u \right\},\nonumber \\
\mathscr{V} &:= \left\{\vm{v}(\vm{x})\in
[ C^0(\Omega) ]^d : \vm{v} \in [ \mathcal{W}(\Omega)]^d \subseteq
[ H^{1}(\Omega) ]^d, \ \vm{v} = \vm{0} \ \textrm{on } \Gamma_u
\right\},\nonumber
\end{align}
\end{subequations}
where the space $\mathcal{W}(\Omega)$ includes linear displacement
fields. The domain is partitioned into elements $\Omega^h$, and on using shape functions $\phi_a$ that span at least the linear space, we substitute vector-valued trial and test functions $\vm{u}^h=\sum_a N_a \vm{u}_a$ and $\vm{v}^h=\sum_b N_b \vm{v}_b$, respectively, into~\eref{eq18} and apply a standard Galerkin procedure to obtain the discrete weak form: find $\vm{u}^h \in \mathscr{U}^h$ such that
\begin{equation}\label{eq:discweakform}
\forall \vm{v}^h \in \mathscr{V}^h \quad a(\vm{u}^h,\vm{v}^h)= \ell(\vm{v}^h),
\end{equation}
which leads to the following system of linear equations:
\begin{subequations}\label{eq:weakform_disc}
\begin{align}
\mat{K}\mat{u}^h&=\mat{f},\\
\mat{K}=\sum_h \mat{K}^h &=\sum_h\int_{\Omega^h}\mat{B}^\transpose\mat{C}\mat{B}\,\rmd \Omega,\\
\mat{f} = \sum_h \mat{f}^h &=\sum_h\left(\int_{\Omega^h}\mat{N}^\transpose\vm{b}\,\rmd \Omega + \int_{\Gamma_t^h}\mat{N}^\transpose\hat{\vm{t}}\,\rmd \Gamma \right),
\end{align}
\end{subequations}
where $\mat{K}$ is the assembled stiffness matrix, $\mat{f}$ the assembled nodal force vector, $\mat{u}^h$ the assembled vector of nodal displacements, $\mat{N}$ is the matrix of shape functions,
$\mat{C}$ is the constitutive matrix for an isotropic linear elastic
material, and $\mat{B}=\bsym{\nabla}_s\mat{N}$ is the strain-displacement matrix that is computed using the derivatives of the shape functions.

\subsection{eXtended Finite Element Method}
\begin{figure}
\centering
\includegraphics[width=0.7\textwidth]{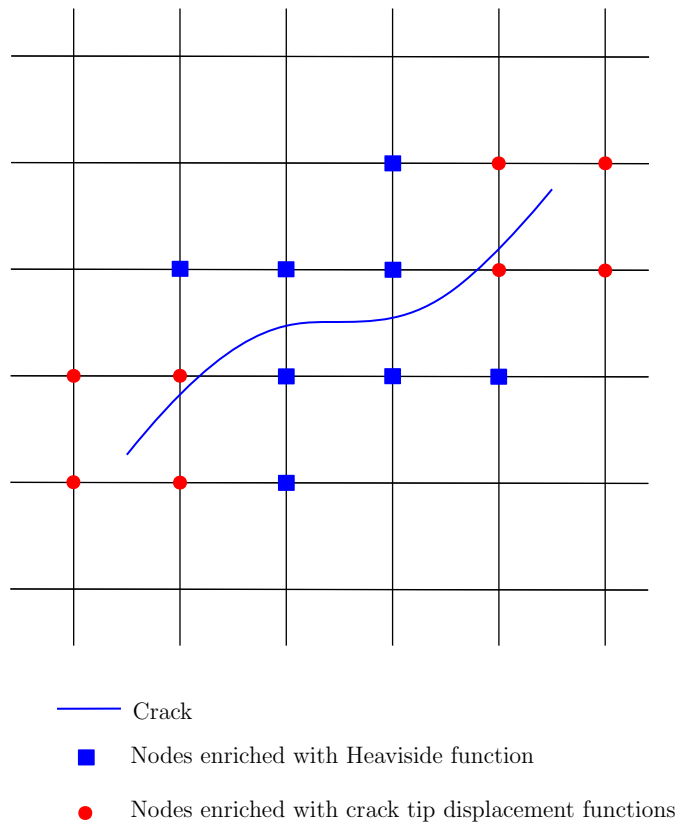}
\caption{XFEM discretisation of a domain with internal discontinuity.}
\label{xfem_discrete}
\end{figure}

With the regular FEM, the mesh topology has to conform to the discontinuous surface. The introduction of the XFEM has alleviated these difficulties by allowing the discontinuities to be independent of the underlying discretization. Within the framework of the eXtended Finite Element Method (XFEM), the trial functions take the following form:
\begin{equation}
\mat{u}^h(\mathbf{x})=\sum_{I\in \mathscr{N}^{std}}\limits \mathbf{N}_I(\xx) \mathbf{u}_I+
  	   \sum_{J\in \mathscr{N}^{hev}}\limits \mathbf{N}_J(\xx) H(\mathbf{x}) \mathbf{a}_J+
       \sum_{K\in \mathscr{N}^{tip}}\limits \mathbf{N}_K(\xx) \left( \sum_{m=1}^4 \limits F_m(r,\theta) \mathbf{b}_K^m\right) 
\label{eq26}
\end{equation}
where $I$ is the set of all the nodes in the system, $J$ is the set of nodes which are completely cut by the crack, $K$ is the set of nodes which contain the crack tips as shown in Figure \ref{xfem_discrete}. $\mathbf{N}_I(\mathbf{x})$ are the standard shape functions associated with the standard DOF $\mathbf{u}_I$, $H(\mathbf{x})$ is the Heaviside function associated with the enriched DOF, $\mathbf{a}_J$ and $F_m(r,\theta)$ are the tip enrichment functions associated with the DOF, $\mathbf{b}_K^m$ that span the near tip asymptotic fields:
\begin{equation}
\begin{aligned}
F_m(r,\theta)&=
\begin{Bmatrix}
\sqrt{r}\sin\dfrac{\theta}{2} & \sqrt{r}\cos\dfrac{\theta}{2} &
\sqrt{r}\sin\theta \sin\dfrac{\theta}{2} & \sqrt{r}\sin\theta \cos\dfrac{\theta}{2}
\end{Bmatrix}\\
\end{aligned}
\label{eq24}
\end{equation}
Following the Galerkin procedure, this modification to the trial function space leads to an enlarged problem to solve:
\begin{equation}
\mathbf{K}\mat{u}^h =\mathbf{F}
\label{eq28}
\end{equation}
where
\begin{equation}
\mathbf{K}^e=
\begin{bmatrix} 
\mathbf{K}^{uu} & \mathbf{K}^{ua} & \mathbf{K}^{ub}\\
\mathbf{K}^{au} & \mathbf{K}^{aa} & \mathbf{K}^{ab}\\
\mathbf{K}^{bu} & \mathbf{K}^{ba} & \mathbf{K}^{bb}
\end{bmatrix}
\label{eq29}
\end{equation}
where, the superscript $uu$ refers to standard FEM components, $aa$ refers to the Heaviside enrichment terms, $bb$ refers to the asymptotic enrichment terms and other terms can be defined similarly. The augmentation of non-polynomial functions to the trial function space, makes the numerical integration non-trivial. This has been of particular interest among research community, for example, equivalent polynomial approach by Ventura~\cite{ventura2006elimination} and Ventura \textit{et al.,}~\cite{ventura2009fast}, conformal mapping~\cite{natarajan2010integrating}, Duffy transformation~\cite{mousavisukumar2010}, generalized Gaussian quadrature~\cite{mousavisukumar2010a}, strain smoothing technique~\cite{bordas2011performance}, exponentially convergent mapping~\cite{minnebo2012},polar mapping~\cite{parkpereira2009} and very recently by using Euler's homogeneous function theorem and Stoke's theorem~\cite{chinlasserre2016}. In~\cite{bordas2011performance}, the strain smoothing technique was combined with the XFEM, coined as smoothed XFEM (Sm-XFEM) to integrate over elements intersected with discontinuous surface. The main advantages of the Sm-XFEM are that no subdivision of the split elements is required and that the derivatives of the shape functions (including the enrichment functions) are not required. However, it was observed that the error level was greater when compared to the conventional XFEM, whilst yielding similar convergence rates.

\section{Linear smoothing in the XFEM}
\label{ssm}

\begin{figure}[htpb]
\centering
        \subfigure[Sub-cell-1]{%
            \label{triangle-1}
            \includegraphics[width=0.3\textwidth]{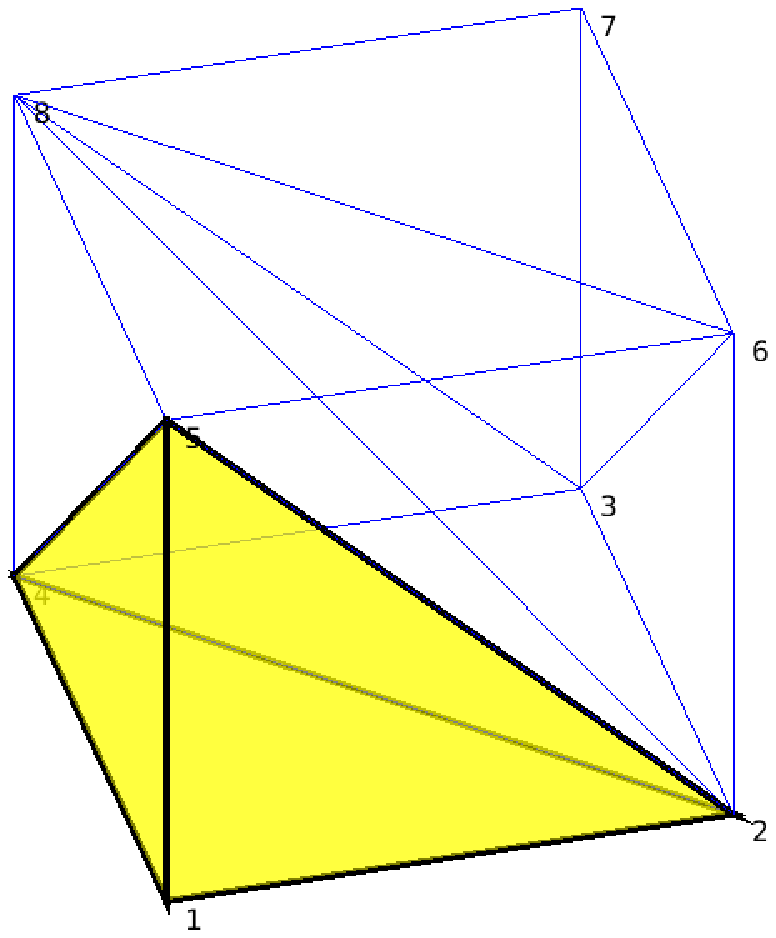}
        }%
        \subfigure[Sub-cell-2]{%
            \label{triangle-2}
            \includegraphics[width=0.3\textwidth]{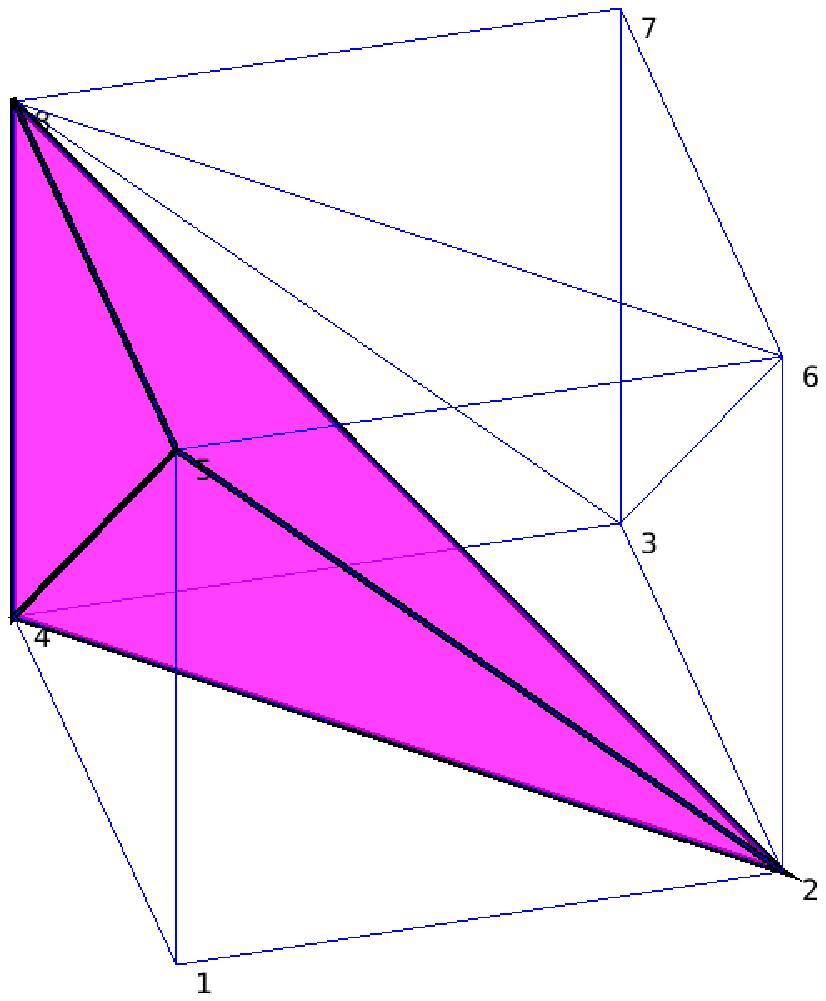}
        }%
        \subfigure[Sub-cell-3]{%
            \label{triangle-3}
            \includegraphics[width=0.3\textwidth]{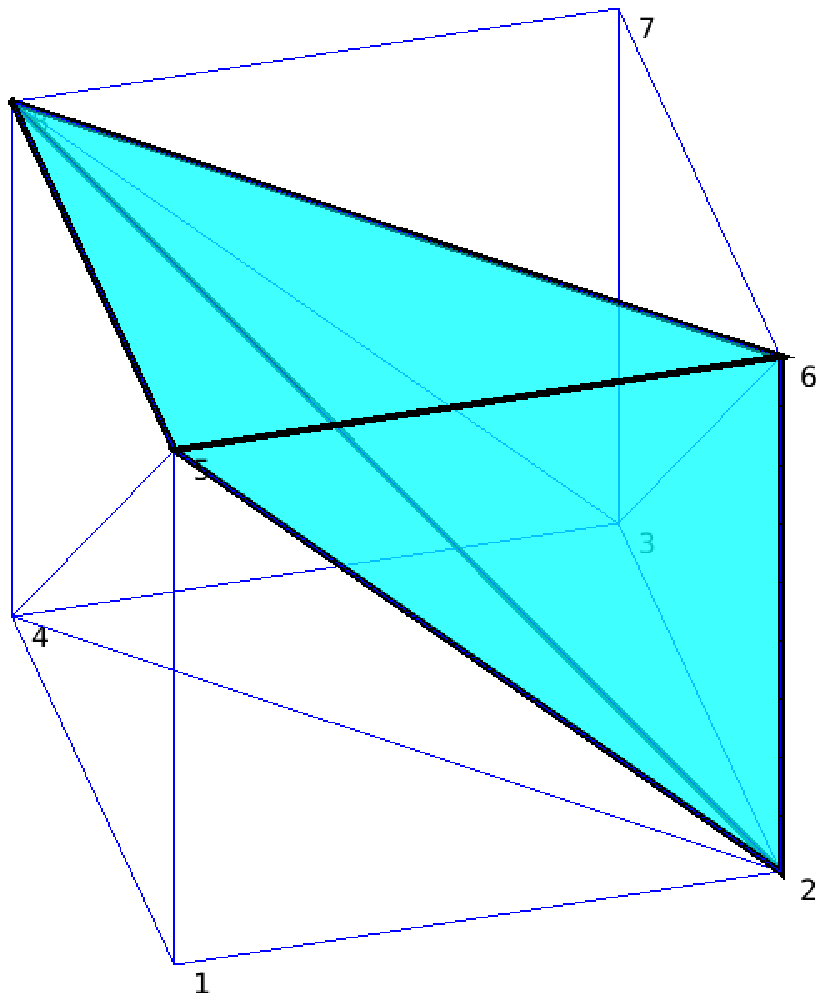}
        }%
        \\ 
        \subfigure[Sub-cell-4]{%
            \label{triangle-4}
            \includegraphics[width=0.3\textwidth]{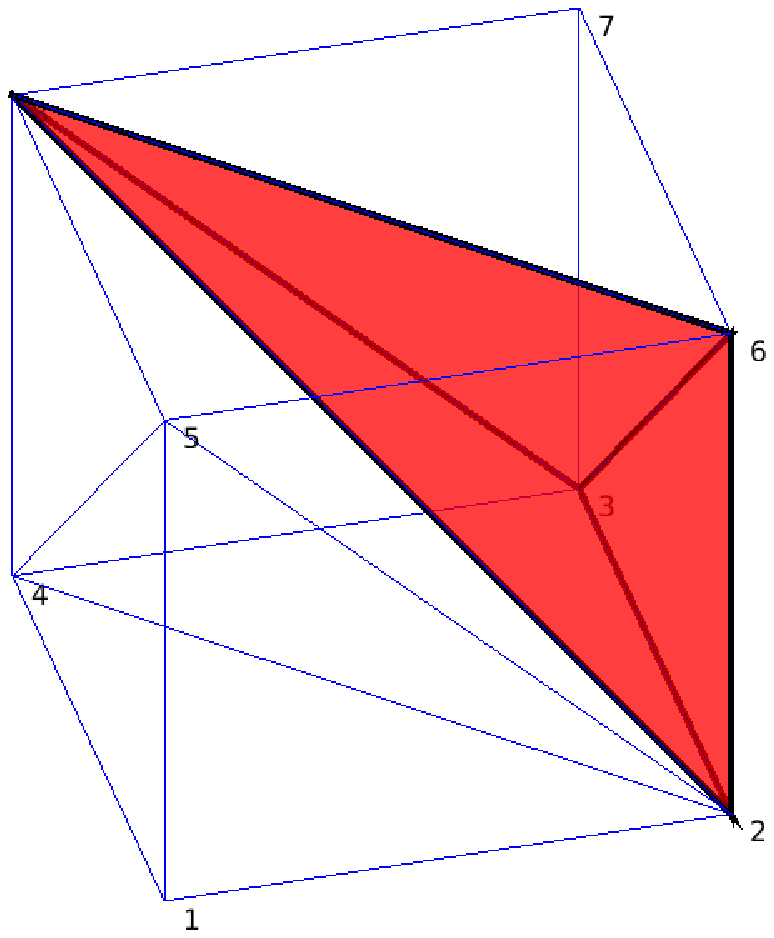}
        }%
        \subfigure[Sub-cell-5]{%
            \label{triangle-5}
            \includegraphics[width=0.3\textwidth]{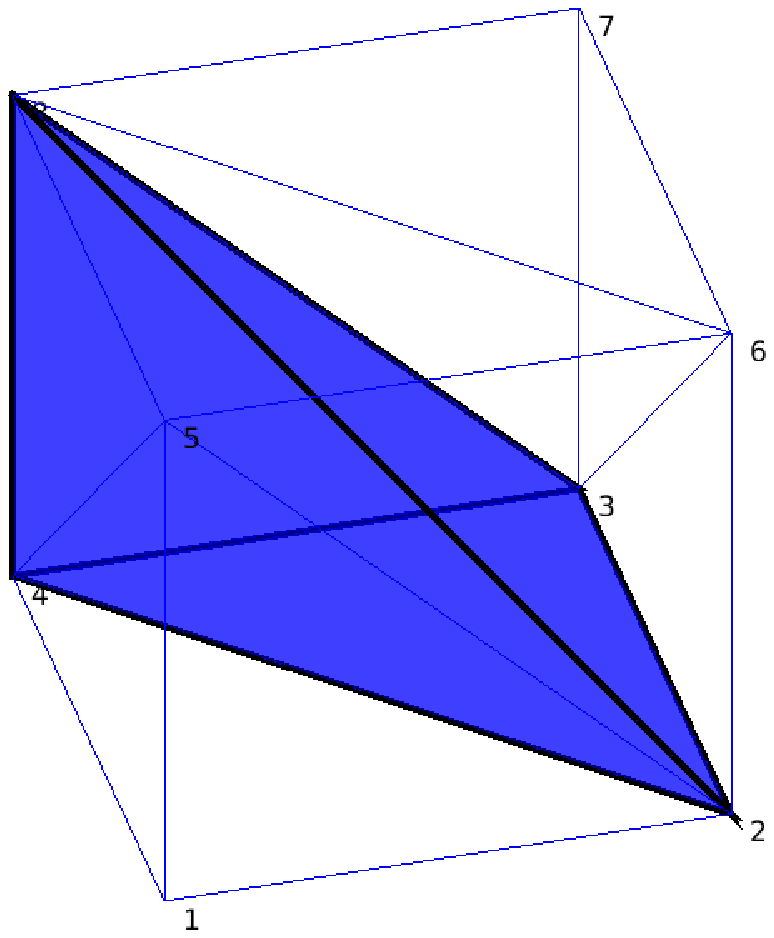}
        }%
        \subfigure[Sub-cell-6]{%
            \label{triangle-6}
            \includegraphics[width=0.3\textwidth]{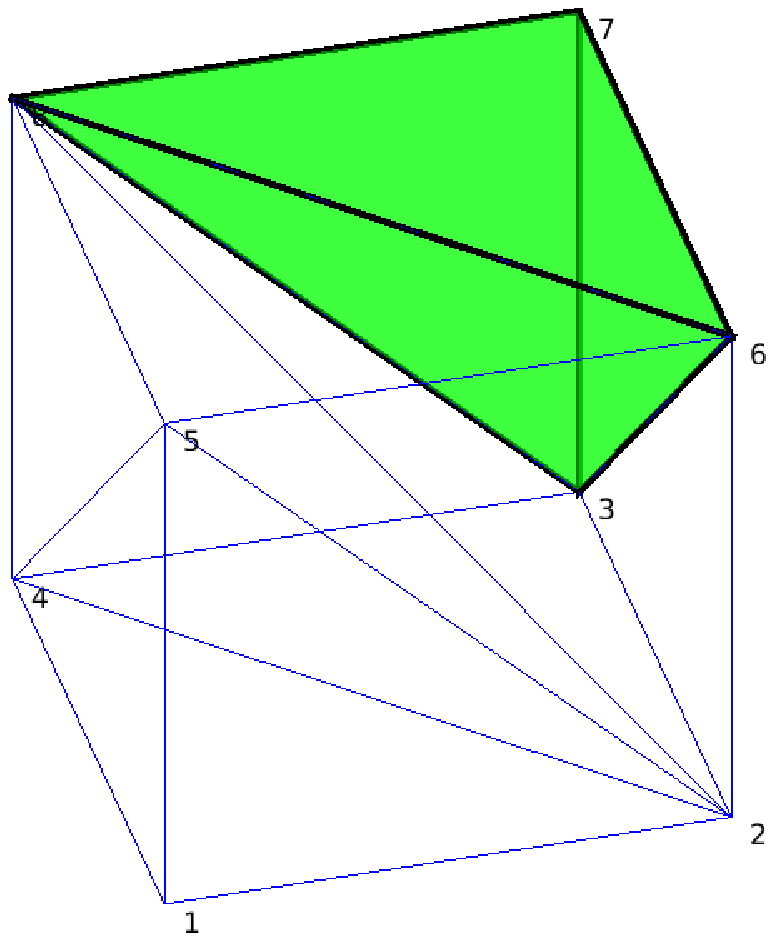}
        }%
\caption{Subdivision of a hexahedral elements into tetrahedral elements. This sub-division is solely for the purpose of numerical integration. A smoothed strain field is computed over each sub-division depending on the choice of smoothing function.}
\label{triangle}
\end{figure}

The strain smoothing was introduced in \cite{chen2001stabilized} for the meshfree methods, which was later extended to the FEM by Liu and co-workers~\cite{liu2007smoothed}. The basic idea is to compute a strain field, referred to as `smoothed' strain field by evaluating the weighted average of the standard (or compatible) strain field. The support domain of the associated material point can be chosen based on surrounding cells, nodes or edges. In this paper, we restrict our discussion only to the cell based strain smoothed FEM. Within this framework, at a point $x_c$ in element $\Omega^h$  the smoothed strain is given below
\begin{equation}
\tilde{\varepsilon}^h_{ij}(\mathbf{x}_c)=
\int_{\Omega^h}\limits{\varepsilon^h_{ij}(\mathbf{x})\Phi(\mathbf{x}-\mathbf{x}_c)\mathrm{d}\Omega}
\label{eq1}
\end{equation}
In terms of the standard element shape function derivatives, $N^h_{I,i}(\mathbf{x})$, the smoothed derivatives are given by
\begin{equation}
\tilde{N}^h_{I,i}(\mathbf{x})=
\int_{\Omega^h}\limits{N^h_{I,i}(\mathbf{x})\Phi(\mathbf{x})\mathrm{d}\Omega}
\label{eq2}
\end{equation}
where, $\Phi$ is the smoothing function and $i=x,y,z$. By invoking the Divergence theorem Eq. 2 can be written as
\begin{equation}
\int_{\Omega^h}\limits{N^h_{I,i}(\mathbf{x})\Phi(\mathbf{x})\mathrm{d}\Omega}=
\int_{\Gamma^h}\limits{N^h_{I}(\mathbf{x})\Phi(x)\mathbf{n(x)}\mathrm{d}\Gamma}
\label{eq3a}
\end{equation}
This form of the strain has the following advantages
\begin{itemize}
\item Domain integration is reduced to a boundary integration along the smoothing cells
\item Does not require the derivatives of the shape functions and hence does not need the Jacobian
\item Does not need isoparametric mapping there by giving a leverage on the distortion level of the mesh
\end{itemize}
The choice of the smoothing function and the integration order used, decide the accuracy of the smoothed strain field. If a constant smoothing function is used, the method is termed the SFEM. It was shown in \cite{bordas2011performance, natarajan2014convergence} that the gradient becomes inaccurate for non-polynomials and higher order polynomial functions. Same issue was also faced with in the context of meshfree approximations and in \cite{duan2012second} this inaccuracy was addressed by introducing an additional domain integral term which ensures consistency between the shape functions and their derivatives. This modified equation was termed the Divergence Consistency (DC). It was also shown that such consistency requirement is implicitly satisfied, if linear field is used. It can be seen that Equation\eqref{eq3b} would reduce to Equation\eqref{eq3a}, if $\Phi$ is a constant.
\begin{equation}
\int_{\Omega^h}\limits{N^h_{I,i}(\mathbf{x})\Phi(\mathbf{x})\mathrm{d}\Omega}=
\int_{\Gamma^h}\limits{N^h_{I}(\mathbf{x})\Phi(x)\mathbf{n(x)}\mathrm{d}\Gamma}-
\int_{\Omega^h}\limits{N^h_{I}(\mathbf{x})\Phi'(\mathbf{x})\mathrm{d}\Omega}
\label{eq3b}
\end{equation}
where, $\Gamma^h$ is the contour of the smoothing cell. Here the domain integral term vanishes as the smoothing function is constant over the domain. Since we assumed linear displacement functions, the strain would be a constant and a unique value can be computed using a single equation. Hence requiring just one interior Gau\ss~points. This can be written as
\begin{subequations}
\begin{equation}
\tilde{N}^h_{I,i}(\mathbf{x}_c)=
\dfrac{1}{A_c}\int_{\Gamma^h}\limits{N^h_{I}(\mathbf{x})\mathbf{n(x)}\mathrm{d}\Gamma}
\label{eq4a}
\end{equation}

\begin{equation}
\tilde{\varepsilon}^h(\mathbf{x}_c)=\tilde{B}_c(\mathbf{x}_c)\mathbf{q}
\label{eq4b}
\end{equation}

\begin{equation}
\tilde{B}_c=
\left[\tilde{\mathbf{B}}_{c1} \quad \tilde{\mathbf{B}}_{c2} \quad \cdots \quad \tilde{\mathbf{B}}_{nc}\right] 
\label{eq4c}
\end{equation}

\begin{equation}
\tilde{N}^h_{I,i}(\mathbf{x}_c)=
\dfrac{1}{V_c}\sum_{b=1}^{nb}\limits
{
\begin{pmatrix} 
N^h_{I}(\mathbf{x}_b^G)\mathbf{n_x} & 0 & 0\\
0 & N^h_{I}(\mathbf{x}_b^G)\mathbf{n_y} & 0\\
0 & 0 & N^h_{I}(\mathbf{x}_b^G)\mathbf{n_z}\\
N^h_{I}(\mathbf{x}_b^G)\mathbf{n_y} & N^h_{I}(\mathbf{x}_b^G)\mathbf{n_x} & 0\\
0 & N^h_{I}(\mathbf{x}_b^G)\mathbf{n_z} & N^h_{I}(\mathbf{x}_b^G)\mathbf{n_y}\\
N^h_{I}(\mathbf{x}_b^G)\mathbf{n_z} & 0 & N^h_{I}(\mathbf{x}_b^G)\mathbf{n_x}\\
\end{pmatrix}
A_b^c
}
\label{eq4d}
\end{equation}
\label{eq4}
\end{subequations}
where, $n_c$ is the number of sub-cells in an element, $V_c$ is the volume of the sub-cell, $n_b$ is the number of boundary surfaces of the sub-cell, $A_b^c$  and $x_b^G$ are the area and Gau\ss~ point of the boundary surface $b$. The smoothing technique has been very efficient for polyhedral elements since the polyhedrons can be divided into number of sub-cells (usually tetrahedrons) and the stiffness matrix is summed up over each sub-cell. It can be seen in Equation\eqref{eq4} that the derivatives of the shape functions are not needed in order to evaluate the strains. Hence, the computation of Jacobian is avoided. This also avoids the associated isoparametric mapping. The stiffness matrix is evaluated as in the regular finite element method by replacing the terms in the strain gradient matrix with the terms evaluated above and summing it up over the sub-cells. The constant smoothing technique when applied to elements other than Constant Strain elements (3-noded triangles and 4-noded tetrahedrons) yields results which are bounded by the results of reduced integration procedure (smoothing over one sub-cell) and full integration procedure (smoothing over infinite number of sub-cells). The method is hence not ‘variationally consistent’ for any number of sub-cells other than 1 and $\infty$ \cite{bordas2010approximation}, whereas the linear smoothing procedure is variationally consistent. The constant smoothing and linear smoothing schemes differ in the choice of the smoothing function. In the linear smoothing scheme the basis function used is $f = \left[ 1 \ x \ y \ z \ xy \ yz \ zx \ xyz \right]^T $ in case of hexahedral subcells and $f = \left[ 1 \ x \ y \ z \right]^T $if tetrahedral sub-cells are used. \fref{triangle} shows one possible division of hexahedral elements into tetrahedral elements (also referred as subcells in the literature) for the purpose of numerical integration. The number of terms in the basis function should be consistent with the number of Gau\ss~ points to obtain a unique solution. Since a linear basis function is being used the domain integral term which results as a consequence of the DC does not vanish and hence it has to be computed by using the appropriate order of Gaussian integration. In the case of tetrahedral sub-cells the system of equations for a linear basis would be


\begin{subequations}

\begin{equation}
\int_{\Omega^h}\limits{\mathbf{N}^h_{I,x}(\mathbf{x})\mathrm{d}\Omega}=
\int_{\Gamma^h}\limits{\mathbf{N}^h_{I}(\mathbf{x})n_x\mathrm{d}\Gamma}
\label{eq5a}
\end{equation}

\begin{equation}
\int_{\Omega^h}\limits{\mathbf{N}^h_{I,x}(\mathbf{x})x\mathrm{d}\Omega}=
\int_{\Gamma^h}\limits{\mathbf{N}^h_{I}(\mathbf{x})xn_x\mathrm{d}\Gamma}-
\int_{\Omega^h}\limits{\mathbf{N}^h_{I}(\mathbf{x})\mathrm{d}\Omega}
\label{eq5b}
\end{equation}

\begin{equation}
\int_{\Omega^h}\limits{\mathbf{N}^h_{I,x}(\mathbf{x})y\mathrm{d}\Omega}=
\int_{\Gamma^h}\limits{\mathbf{N}^h_{I}(\mathbf{x})yn_x\mathrm{d}\Gamma}
\label{eq5c}
\end{equation}

\begin{equation}
\int_{\Omega^h}\limits{\mathbf{N}^h_{I,x}(\mathbf{x})z\mathrm{d}\Omega}=
\int_{\Gamma^h}\limits{\mathbf{N}^h_{I}(\mathbf{x})zn_x\mathrm{d}\Gamma}
\label{eq5d}
\end{equation}

\end{subequations}
for $\mathbf{N}^h_{I,x}(\mathbf{x})$ 
\begin{subequations}
\begin{equation}
\int_{\Omega^h}\limits{\mathbf{N}^h_{I,y}(\mathbf{x})\mathrm{d}\Omega}=
\int_{\Gamma^h}\limits{\mathbf{N}^h_{I}(\mathbf{x})n_y\mathrm{d}\Gamma}
\label{eq6a}
\end{equation}

\begin{equation}
\int_{\Omega^h}\limits{\mathbf{N}^h_{I,y}(\mathbf{x})x\mathrm{d}\Omega}=
\int_{\Gamma^h}\limits{\mathbf{N}^h_{I}(\mathbf{x})xn_y\mathrm{d}\Gamma}
\label{eq6b}
\end{equation}

\begin{equation}
\int_{\Omega^h}\limits{\mathbf{N}^h_{I,y}(\mathbf{x})y\mathrm{d}\Omega}=
\int_{\Gamma^h}\limits{\mathbf{N}^h_{I}(\mathbf{x})yn_y\mathrm{d}\Gamma}-
\int_{\Omega^h}\limits{\mathbf{N}^h_{I}(\mathbf{x})\mathrm{d}\Omega}
\label{eq6c}
\end{equation}

\begin{equation}
\int_{\Omega^h}\limits{\mathbf{N}^h_{I,y}(\mathbf{x})z\mathrm{d}\Omega}=
\int_{\Gamma^h}\limits{\mathbf{N}^h_{I}(\mathbf{x})zn_y\mathrm{d}\Gamma}
\label{eq6d}
\end{equation}

\end{subequations}

for $\mathbf{N}^h_{I,y}(\mathbf{x})$. 


\begin{subequations}

\begin{equation}
\int_{\Omega^h}\limits{\mathbf{N}^h_{I,z}(\mathbf{x})\mathrm{d}\Omega}=
\int_{\Gamma^h}\limits{\mathbf{N}^h_{I}(\mathbf{x})n_z\mathrm{d}\Gamma}
\label{eq6az}
\end{equation}

\begin{equation}
\int_{\Omega^h}\limits{\mathbf{N}^h_{I,z}(\mathbf{x})x\mathrm{d}\Omega}=
\int_{\Gamma^h}\limits{\mathbf{N}^h_{I}(\mathbf{x})xn_z\mathrm{d}\Gamma}
\label{eq6bz}
\end{equation}

\begin{equation}
\int_{\Omega^h}\limits{\mathbf{N}^h_{I,z}(\mathbf{x})y\mathrm{d}\Omega}=
\int_{\Gamma^h}\limits{\mathbf{N}^h_{I}(\mathbf{x})yn_z\mathrm{d}\Gamma}
\label{eq6cz}
\end{equation}

\begin{equation}
\int_{\Omega^h}\limits{\mathbf{N}^h_{I,z}(\mathbf{x})z\mathrm{d}\Omega}=
\int_{\Gamma^h}\limits{\mathbf{N}^h_{I}(\mathbf{x})zn_z\mathrm{d}\Gamma}-
\int_{\Omega^h}\limits{\mathbf{N}^h_{I}(\mathbf{x})\mathrm{d}\Omega}
\label{eq6dz}
\end{equation}
\end{subequations}

\begin{figure}[htpb]
\centering
\includegraphics[width=0.9\textwidth]{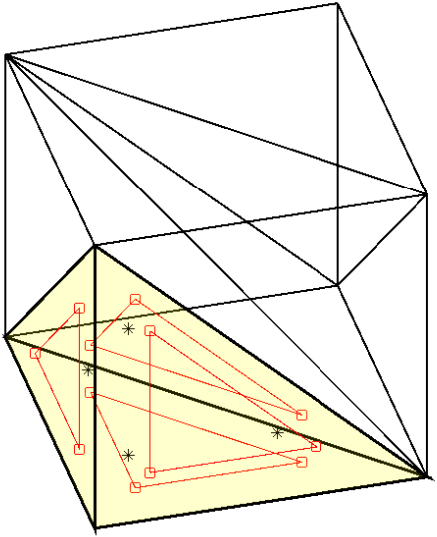}
\caption{The location of Gau\ss~points for boundary integration and domain integration over a tetrahedron sub-cell of a hexahedral element.}
\label{gps}
\end{figure}
for $\mathbf{N}^h_{I,z}(\mathbf{x})$. Here $\mathbf{N}_I$ represents the shape function associated with the $I^{th}$ node of the parent element. It is independent of the sub-cell. The location of the gauss points for the boundary integration and domain integration in a tetrahedral sub-cell are shown in  Figure \ref{gps}. Let the natural coordinates of the $m^{th}$ interior gauss point of a sub-cell be $\mathbf{p}_m=(x_m, y_m, z_m)$ and its associated gauss weight be $w_m$; coordinates of the $k^{th}$ boundary of the sub-cell be $\mathbf{c}^k_g=(x^k_g, y^k_g, z^k_g)$ and the associated weights be $v^k_g$ . The unit outward normal associated with the $g^{th}$ gauss point of the $k^{th}$ boundary of the sub-cell is denoted by $\mathbf{n}^k=(n^k_x,n^k_y,n^k_z )$. The smoothed derivatives are computed numerically as follows

\begin{equation}
\mathbf{W}\mathbf{d}_i=\mathbf{f}_i \quad
\textrm{where,} \ i=x,y,z
\label{eqn:eq7}
\end{equation}

\begin{equation}
\mathbf{W}=
\begin{pmatrix} 
w_1 & w_2 & w_3 & w_4\\
w_1x_1 & w_2x_2 & w_3x_3 & w_4x_4\\
w_1y_1 & w_2y_2 & w_3y_3 & w_4y_4\\
w_1z_1 & w_2z_2 & w_3z_3 & w_4z_4\\
\end{pmatrix}
\label{eqn:eq8}
\end{equation}

\begin{equation}
\mathbf{f}_x=
\begin{pmatrix*}[l] 
\sum_{k=1}^4\limits \sum_{g=1}^3\limits N_I(\mathbf{c}_g^k) n_x^k v_g^k\\
\sum_{k=1}^4\limits \sum_{g=1}^3\limits N_I(\mathbf{c}_g^k) x_g^k n_x^k v_g^k -
\sum_{m=1}^4\limits N_I(\mathbf{p}_m)w_m\\
\sum_{k=1}^4\limits \sum_{g=1}^3\limits N_I(\mathbf{c}_g^k) y_g^k n_x^k v_g^k\\
\sum_{k=1}^4\limits \sum_{g=1}^3\limits N_I(\mathbf{c}_g^k) z_g^k n_x^k v_g^k
\end{pmatrix*}
\label{eqn:eq9}
\end{equation}

\begin{equation}
\mathbf{f}_y=
\begin{pmatrix*}[l] 
\sum_{k=1}^4\limits \sum_{g=1}^3\limits N_I(\mathbf{c}_g^k) n_y^k v_g^k\\
\sum_{k=1}^4\limits \sum_{g=1}^3\limits N_I(\mathbf{c}_g^k) x_g^k n_y^k v_g^k\\
\sum_{k=1}^4\limits \sum_{g=1}^3\limits N_I(\mathbf{c}_g^k) y_g^k n_y^k v_g^k-
\sum_{m=1}^4\limits N_I(\mathbf{p}_m) w_m\\
\sum_{k=1}^4\limits \sum_{g=1}^3\limits N_I(\mathbf{c}_g^k) z_g^k n_y^k v_g^k
\end{pmatrix*}
\label{eqn:eq10}
\end{equation}

\begin{equation}
\mathbf{f}_z=
\begin{pmatrix*}[l] 
\sum_{k=1}^4\limits \sum_{g=1}^3\limits N_I(\mathbf{c}_g^k) n_z^k v_g^k\\
\sum_{k=1}^4\limits \sum_{g=1}^3\limits N_I(\mathbf{c}_g^k) x_g^k n_z^k v_g^k\\
\sum_{k=1}^4\limits \sum_{g=1}^3\limits N_I(\mathbf{c}_g^k) y_g^k n_z^k v_g^k\\
\sum_{k=1}^4\limits \sum_{g=1}^3\limits N_I(\mathbf{c}_g^k) z_g^k n_z^k v_g^k-
\sum_{m=1}^4\limits N_I(\mathbf{p}_m) w_m\\
\end{pmatrix*}
\label{eqn:eq10a}
\end{equation}

The smoothed derivative of the $i^{th}$ shape function evaluated at the four interior Gau\ss~points of a tetrahedral sub-cell is given by

\begin{equation}
\begin{aligned}
d_i&=
\begin{bmatrix}
d_i^1 & d_i^2 & d_i^3 & d_i^4
\end{bmatrix}^T \\
&=
\begin{bmatrix}
N_{I,i}(p_1) & N_{I,i}(p_2) & N_{I,i}(p_3) & N_{I,i}(p_4)
\end{bmatrix}^T \textrm{where,} \ i=x,y,z\\
\end{aligned}
\label{eqn:eq11}
\end{equation}

The same procedure is to be repeated for all the shape functions of the parent element. For the $m^{th}$ interior gauss point of a sub-cell of a $n$ sided polygon the smoothed strain displacement matrix is given by

\begin{equation}
\begin{aligned}
\tilde{\mathbf{B}}_c(\mathbf{p}_m)&=
\begin{bmatrix}
\tilde{\mathbf{B}}_{c1} (\mathbf{p}_m) & \tilde{\mathbf{B}}_{c2} (\mathbf{p}_m) & \cdots & \tilde{\mathbf{B}}_{c3}(\mathbf{p}_m)
\end{bmatrix} & \textrm{where,} &\ m=1,2,3,4\\
\end{aligned}
\label{eqn:eq12}
\end{equation}

\begin{equation}
\tilde{\mathbf{B}}_{cI}(p_m) =
\begin{bmatrix*}[l]
d_x^1 & 0 & 0\\
0 & d_y^1 & 0\\
0 & 0 & d_z^1\\ 
d_y^1 & d_x^1 & 0\\
0 & d_z^1 & d_y^1\\
d_z^1 & 0 & d_x^1\\
\end{bmatrix*}
\label{eqn:eq13}
\end{equation}

For the displacement approximation given by \eref{eq26}, the compatible strain field is given by:
\begin{equation}
\varepsilon^h(\mathbf{x})= 
\begin{bmatrix}
\mathbf{B}_{\rm fem} &  \mathbf{B}_{\rm hev} & \mathbf{B}_{\rm tip}
\end{bmatrix} \mathbf{q}^T
\label{eq31b}
\end{equation}
where $\mathbf{q} = \lbrace \mathbf{u} \ \mathbf{a} \ \mathbf{b} \rbrace$ is the vector of degrees of freedom, $\mathbf{B}_{\rm fem}$, $\mathbf{B}_{\rm hev}$ and $\mathbf{B}_{\rm tip}$  contains the strain displacement terms corresponding to the regular finite element part, Heaviside enriched part and the tip enriched part. The components of the compatible strain field are:
\begin{align}
\mathbf{B}_{\rm fem} &= \mathbf{L} \mathbf{N}_I \nonumber \\
\mathbf{B}_{\rm hex} &= \mathbf{L}\mathbf{N}_J \left( H(\mathbf{x})-H(\mathbf{x}_J) \right) \nonumber \\
\mathbf{B}_{\rm tip} &= \mathbf{L} \mathbf{N}_K \left( \sum\limits_{m=1}^4 \left( F_m(\mathbf{x}) - F_m(\mathbf{x}_K) \right) \right)
\end{align}
where,
\begin{equation*}
\mathbf{L} = \begin{bmatrix*}[l]
\frac{\partial}{\partial x} & 0 & 0 \\
 0 & \frac{\partial}{\partial y} & 0 \\
 0 & 0 & \frac{\partial}{\partial z}\\
\frac{\partial}{\partial y} & \frac{\partial }{\partial x} & 0\\
0 & \dfrac{\partial}{\partial z} & \frac{\partial } {\partial y}\\
\frac{\partial }{\partial z} & 0 & \frac{\partial}{\partial x}\\
\end{bmatrix*}
\end{equation*}
The smoothed counterpart of the above compatible strain field can be obtained by following the procedure outlined earlier. The elements that are intersected by the discontinuous surface is divided into tetrahedra and a linear smoothing basis, $f(\mathbf{x}) = \left\lbrace 1 \ x \ y \ z \right\rbrace $  is chosen to evaluate the smoothed strain.

\begin{remark}
In case of two dimensions, the subcell is a triangle and the smoothing procedure can be derived from the linear basis
\begin{equation}
f(\xx) = \left\lbrace 1 \ x \ y \right\rbrace
\end{equation}
with derivative
\begin{equation}
f_{,j}(\xx) = \left\lbrace 0 \ \delta_{1j} \ \delta_{2j} \right\rbrace
\end{equation}
\end{remark}


\section{Numerical Examples}
\label{numerics}

In this section, the accuracy and the convergence properties of the proposed formulation is numerically studied within the framework of linear elastic fracture mechanics (LEFM) in both two and three dimensions. The domain is discretized with four noded quadrilateral and eight noded hexahedral elements in two and three dimensions, respectively. The numerical results from the present formulation is compared with the conventional XFEM and the SmXFEM~\cite{bordas2011performance}. The following convention is adopted to compute the stiffness matrix within the framework of the smoothing technique:

\begin{table}[htpb]
\caption{Number of sub-cells used to compute the stiffness matrix for the constant smoothed XFEM (Sm-XFEM) and the linear smoothed XFEM (LSm-XFEM). In case of Sm-XFEM, the smoothing function is chosen as $f(\mathbf{x}) = 1$, whilst in case of LSm-XFEM, a complete set of polynomials is chosen. For example, $f(\mathbf{x}) = \lbrace 1 \ x \ y \rbrace$ for two dimensions and $f(\mathbf{x}) = \lbrace 1 \ x \ y \ z \rbrace$ for three dimensions as smoothing function.}  
\centering
\begin{tabular}{clrr} 				
\hline		
& \\					
& Type of element & Sm-XFEM & LSm-XFEM \\ 
& \\
\hline \\[0.1ex]
\multirow{3}{*}{two dimensions} & Standard elements & 4 & 1 \\  [0.5ex]
&Tip enriched elements & 5 & 5 \\ [0.5ex]
&Split enriched elements & 8 & 8 \\ [0.5ex]
\cline{2-4} \\[0.1ex]
\multirow{3}{*}{three dimensions} & Standard  elements & 6 & 1 \\ [0.5ex]
&Tip enriched elements & 24 & 24 \\[0.5ex]
& Split enriched elements & 12 & 12 \\ [0.5ex]
\hline
\end{tabular}
\label{table_gauss} 
\end{table}
For the conventional XFEM, the elements that are intersected by the discontinuous surface is triangulated and a higher order triangular quadrature scheme is adopted. And for the standard elements, 2$\times$2 Gau\ss~ quadrature rule is adopted. To estimate the error and to study the convergence properties, the $L^2$ norm and the $H^1$ semi-norm is used.

\subsection{Infinite plate with center crack under far-field tension}
\begin{figure}[htpb]
\centering
\includegraphics[width=0.5\textwidth]{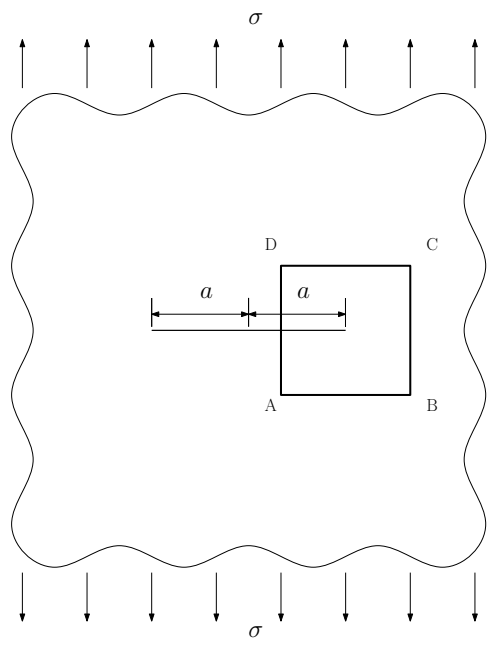}
\caption{An infinite plate with a center crack subject to far-field tensile stress}
\label{fig_grif}
\end{figure}

\begin{figure}[htpb]
\centering
\subfigure[Relative error in $L^2$ norm]{\includegraphics[width=0.9\textwidth]{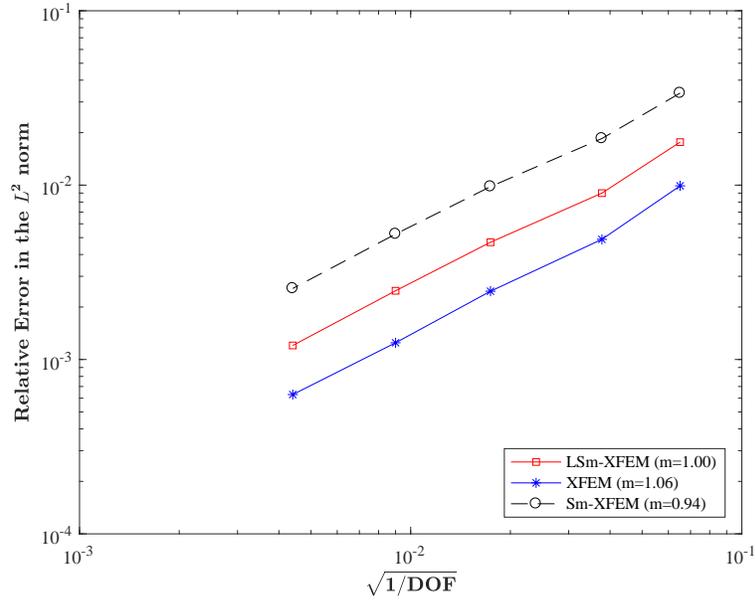}}
\subfigure[Relative error in SIF]{\includegraphics[width=0.9\textwidth]{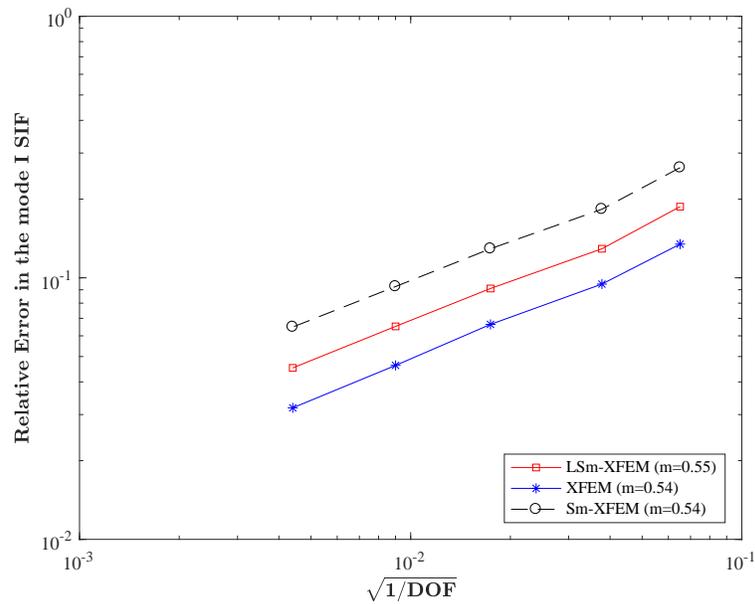}}
\caption{Convergence of the relative error in the displacement and in the stress intensity factor with mesh refinement for a infinite plate with a center crack subjected to uniform tensile stress. On the boundary, Westergaard solution is applied.}
\label{fig_grif_l2_k}
\end{figure}

Consider an infinite plate with a centre crack subjected to far field tension under plane strain condition has been considered. Consider an infinite plate as shown in Figure \ref{fig_grif} A small section ABCD has been solved. The effect of the far-field stress has been accounted by prescribing equivalent displacements as given by following closed form solution Equation \eqref{eq44} in polar coordinates centered at the crack tip.
\begin{align}
u_x(r,\theta) &=\frac{2(1+\nu)}{\sqrt{2\pi}}
		      \frac{K_I}{E}
		      \sqrt{r}\cos\frac{\theta}{2}
		      \left( 2-2\nu-\cos^2 \frac{\theta}{2}\right) \nonumber \\
u_y(r,\theta) &=\frac{2(1+\nu)}{\sqrt{2\pi}}
		      \frac{K_I}{E}
		      \sqrt{r}\sin\frac{\theta}{2}
		      \left( 2-2\nu-\cos^2 \frac{\theta}{2}\right)
\label{eq44} 
\end{align}
where  $K_I$ $=$ $\sigma$ $\sqrt{\pi a}$, the mode I stress intensity factor, $\nu$ is the Poisson's ratio, $E$ is the Young’s modulus. The dimension has been chosen as 10 x 10 mm. ‘a’ is chosen as 100 mm. The convergence of the relative error in the displacement ($L^2$ norm) and the stress intensity factor is shown in Figure \ref{fig_grif_l2_k}. It can be seen that in general all the three methods yields a rate of convergence of 1 in the $L^2$ norm and 0.5 in the $H^1$ semi-norm. For a given dof, the conventional XFEM yields slightly accurate results than the Sm-XFEM or the LSm-XFEM but the errors are within the same order. Moreover it is noted that in the XFEM, 13 integration points per triangle (for the tip element) is used when compared to three integration points in case of LSm-XFEM and one integration point in case of Sm-XFEM. The sub-optimal rate of convergence is due to the fact that we are employing topological enrichment scheme as opposed to geometric enrichment. 

\subsection{Finite Plate with an edge crack subject to tensile and shear stresses}
\begin{figure}[htpb]
\centering
\includegraphics[width=1\textwidth]{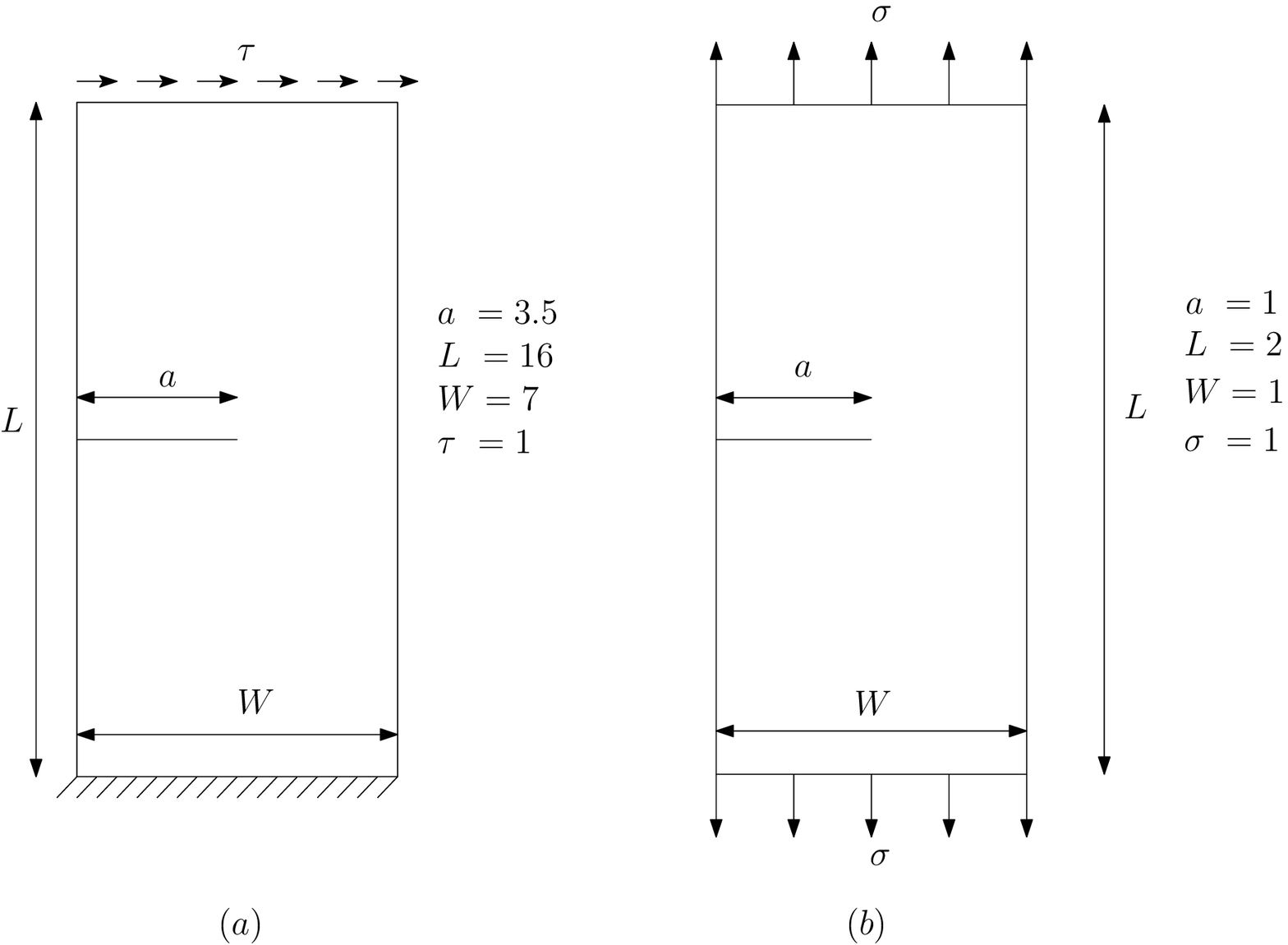}
\caption{Geometry and boundary conditions for a finite plate with an edge crack subject to: (a) uniform shear stress at the top face and (b) uniform tensile stress.}
\label{fig_edge_shear}
\end{figure}

Next, consider a finite element with an edge crack subjected to tensile and shear stresses as shown in Figure \ref{fig_edge_shear}. A consistent system of units is used for the analysis.

\begin{figure}[htpb]
\centering
\includegraphics[width=1\textwidth]{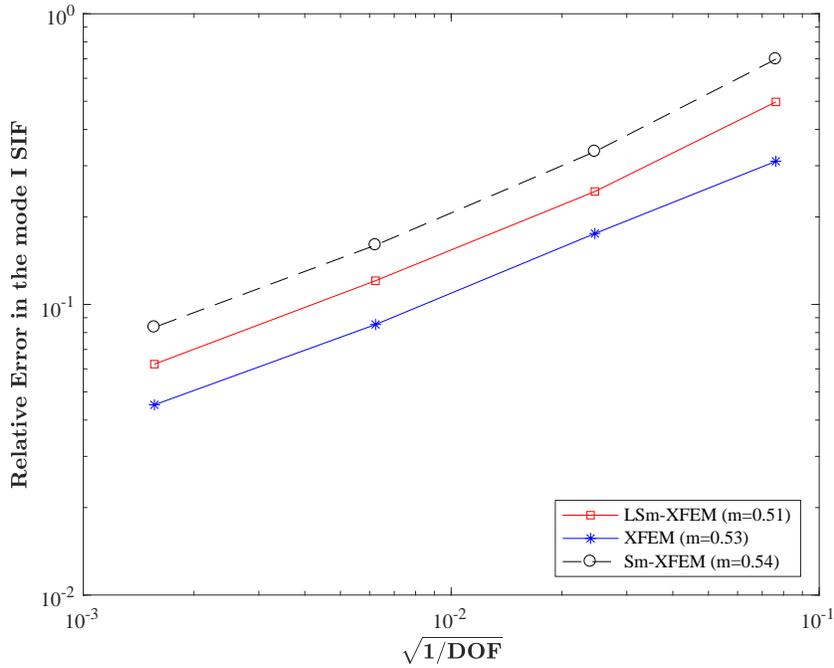}
\caption{Relative error in the mode I SIF for edge cracked plate with tensile loading. }
\label{fig_edge_k1}
\end{figure}

\paragraph{Plate subjected to tensile stress}
In this case, the dimension of the plate is 1 x 2 units. The Youngs' modulus, E and Poisson’s ratio, $ \nu $ are taken as 1000 and 0.3 respectively. A state of plane strain condition is assumed. The crack width is taken as 0.5 units. The obtained SIF are compared with the reference empirical solution\cite{ewalds1984fracture}:
\begin{equation}
K_{ref}=f(\alpha) \ \sigma \sqrt{\pi a}
\label{eq45}
\end{equation}
where, $f(\alpha)=1.12-0.231\alpha+10.55\alpha^2-21.72\alpha^3+30.39\alpha^4$, $ \alpha =a/W $  is the crack width ratio, $a$ is the half-crack width and $w$ is the plate width. The convergence of the relative error in the stress intensity factor is shown in Figure \ref{fig_edge_k1}. It can be seen that the all the three methods converge at almost identical rates ($\approx$ 0.5). The results of LS scheme are better than the CS scheme and are almost equal with the conventional XFEM.

\begin{figure}[htpb]
\centering
\subfigure[]{\includegraphics[width=0.9\textwidth]{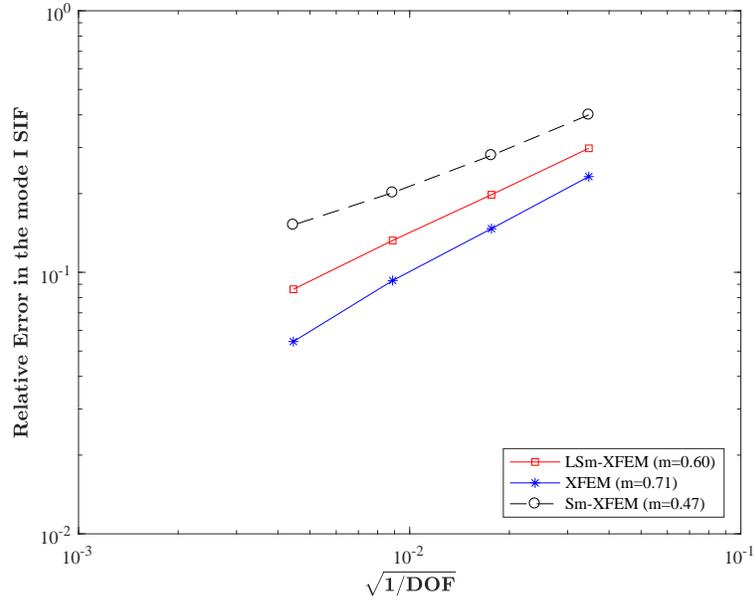}}
\subfigure[]{\includegraphics[width=0.9\textwidth]{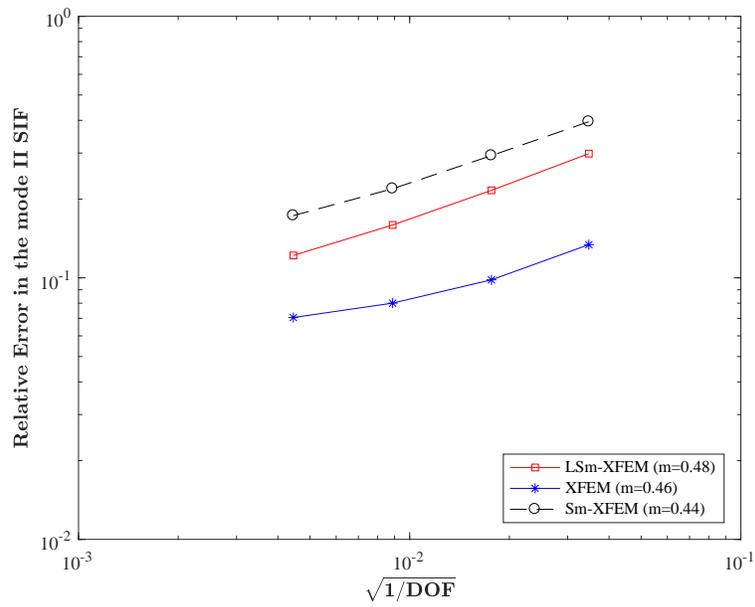}}
\caption{Relative error in the mode I and mode II stress intensity factors for a plate with an edge crack subjected to shear stress.}
\label{shear_k1_k2}
\end{figure}

\paragraph{Plate subjected to shear stresses}
In this case, the dimensions of the plate are taken as $W=$ 7 units and $L=$ 16 units. The plate is subjected to shear stress on the top edges, while the displacements are constrained at the bottom edge. The crack width is taken as 3.5 units. The Youngs' modulus, $E$ and Poisson’s ratio, $\nu$ are taken as 3$ \times$ 10$^7$ and 0.25 respectively. Plane strain condition is assumed. The reference SIF is taken from \cite{ewalds1984fracture}, which is $K_I$ = 34 units, $K_{II}$ = 4.55 units. The convergence of the $K_I$ and the $K_{II}$ are presented in Figure \ref{shear_k1_k2}. It is again seen that all the three methods have similar convergence rates. The LS scheme is also more accurate than the CS scheme with a very minor additional computational expense. It is again recalled that the additional integration points still require only the shape function values which can be obtained by linear interpolation along the boundary. The error can be attributed to the inadequate approximation space in the local crack tip region, i.e, the asymptotic fields are approximated by a linear field.

\subsection{Plate with an inclined center crack subject to tensile stresses}
\begin{figure}[htpb]
\centering
\includegraphics[width=0.5\textwidth]{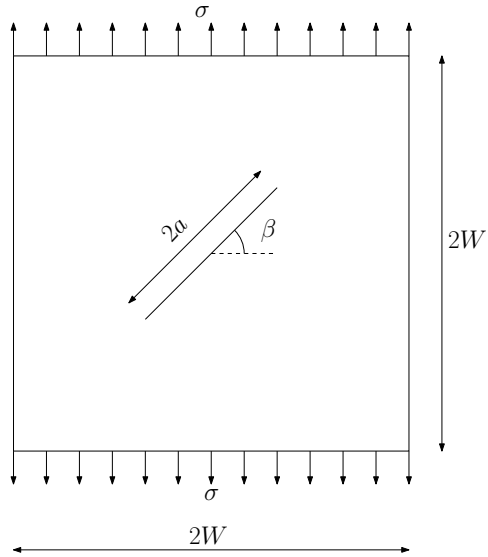}
\caption{Plate with an inclined center crack subject to tensile stress: geometry and boundary conditions}
\label{inclined_plate}
\end{figure}

\begin{figure}
\centering
\includegraphics[width=0.9\textwidth]{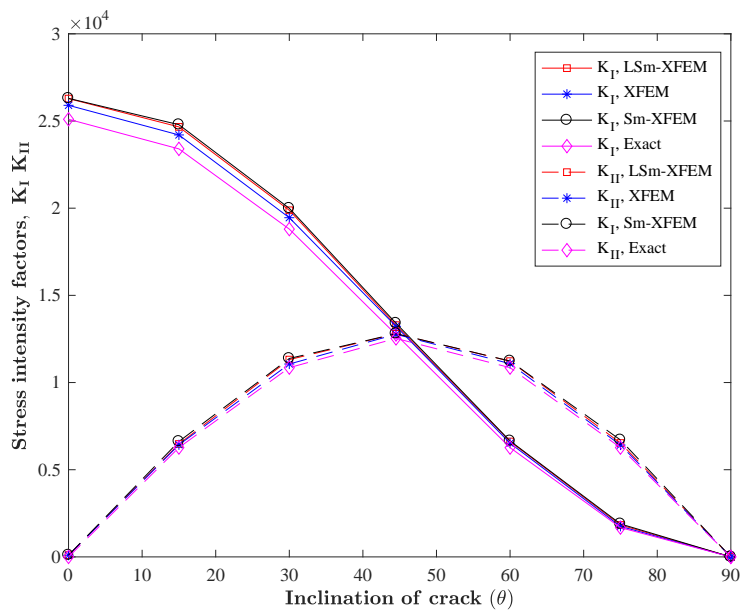}
\caption{Influence of inclination of the crack on the Mode I and mode II stress intensity factors for a plate with a center crack subjected to far field tensile stress.}
\label{inclined_infinite_corr}
\end{figure}

Next, to illustrate the efficacy of the formulation SIFs in case of mixed-mode loading conditions, consider an inclined center crack subjected to far field tension (see \fref{inclined_plate}). The dimensions of the plate are taken as 20 $\times$ 20. The crack width, $2a$ is chosen as 2 units. A far field uniform tensile stress, 1$\times$10$^4$ units is applied with Young's modulus, $E=$ 1$\times$10$^7$ and Poisson's ratio, $\nu=$ 0.3. The accuracy of the numerically computed SIFs are compared with analytical SIFs given by:
\begin{align}
K_{I} &= \sigma \sqrt{\pi a} \cos^2 (\beta) \nonumber \\
K_{II} &= \sigma \sqrt{\pi a} \sin (\beta) \cos (\beta)  
\label{eq48}
\end{align}
where $\beta$ is the inclination of the crack measured anti-clockwise from the $x-$ axis. Based on a progressive refinement, it was observed that a structured mesh of 100 $\times$ 100 quadrilateral mesh is adequate. The influence of the crack angle and different modelling approaches, viz., XFEM, Sm-XFEM, LSm-XFEM on the SIFs are shown in \fref{inclined_infinite_corr}. It can be seen that the results from the proposed approach are accurate and comparable with the conventional XFEM and slightly more accurate than the Sm-XFEM.

\subsection{Finite plate with a through-thickness edge-crack subject to tensile stresses}
\begin{figure}[htpb]
\centering
\includegraphics[width=0.4\textwidth]{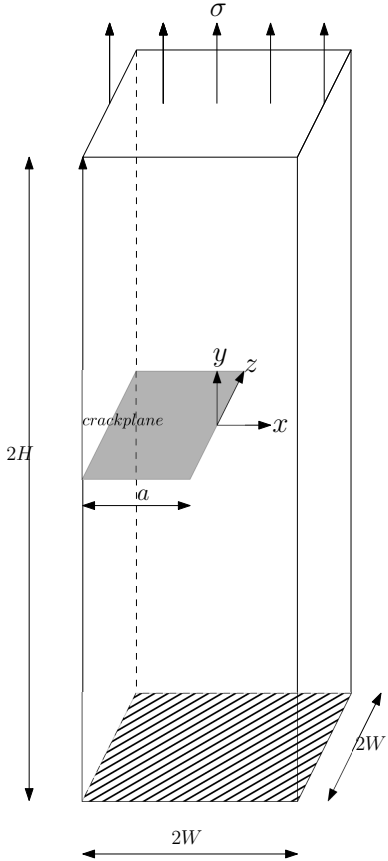}
\caption{Finite plate with a through-thickness edge-crack subject to tensile stress}
\label{prob_3d}
\end{figure}

As a last example, the linear smoothing technique is extended to three dimensional domain with a through-the-thickness edge crack subjected to uniform tensile stress as shown in \fref{prob_3d} with dimensions $W/a=$ 1 and $H/W=$ 3. The displacement at the bottom face is constrained in all directions and a uniform tensile stress $\sigma = $ 1$\times$10$^4$ is applied on the top face. The material properties are: Young's modulus $E=$ 1$\times$10$^7$ and Poisson's ratio $\nu=$ 0.3. The domain is discretized with structured eight noded hexahedral elements and the normalised SIF from \cite{saputra2015computation} is taken as the reference solution. 

\begin{figure}[htpb]
\centering
\includegraphics[width=1\textwidth]{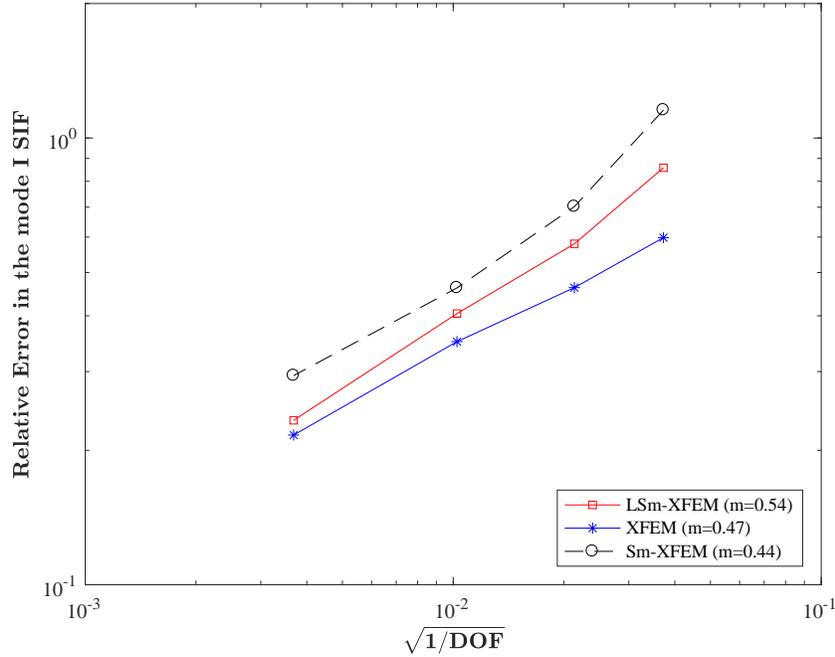}
\caption{Relative error in the normalized SIF  $\dfrac{K_{I}}{\sigma \sqrt{\pi a}}$ for a three-dimensional domain with an edge crack subjected to uniform tensile stress.}
\label{conv_lsm}
\end{figure}

The smoothed strain field over a standard element is computed without any further sub-divisions and with $f(\xx) = \left[ 1 \ x \ y \ z\ xy\ yz \ zx \ xyz \right]$ as a smoothing function. For the elements that are intersected by the discontinuous surface, the element is sub-divided into tetrahedra and $f(\xx) = \left[ 1 \ x \ y \ z\right]$ is chosen as the smoothing function. In case of the LSm-XFEM a total of 96 Gau\ss~ points are used in case of tip enriched elements where as 300 Gau\ss~ points are used in the conventional XFEM. In the case of Heaviside enriched elements 48 Gau\ss~ points are used in case of LSm-XFEM and 60 Gau\ss~ points are used in case of the conventional XFEM. The convergence of the relative error in the normalised stress intensity factor is shown in \fref{conv_lsm}. It can be seen that the LSm-XFEM is more accurate than the Sm-XFEM and is in good agreement with the conventional XFEM.

\section{Conclusions}
\label{conc}
In this paper the Linear smoothing (Second order smoothing) was discussed and a method to couple it with the extended finite element method was presented. The developed method was used to solve problems with discontinuities and singularities in both two and three dimensions. The method also involves a rational integration procedure employing the Green’s theorem. The performance of the linear smoothing scheme for enriched approximation space was studied. Through numerical examples it was shown that the Linear smoothing scheme is more accurate than its constant counterpart. The linear smoothing scheme leads to almost identical results to the standard extended finite element method, but it requires fewer quadrature points, viz., approximately one-fourth to what is required with the conventional XFEM. 

The constant smoothing and the linear smoothing technique is extended to three dimensions for the first time.  Although the presented example in three dimensions is for straight crack, it can be easily extended to other crack profiles. The superior accuracy of the linear smoothing technique is also obtained in the three dimensional case. These results are attributed to the superior approximation properties of the linear smoothing compared to the constant strain smoothing, which is immediately apparent for problems involving complex, non-polynomial, enrichment functions. The remaining, incompressible, error level is attributed to the inadequate approximation space in the smoothed strain, i.e. to the inability of a linear smoothed strain to approximate the singular strains provided by the enriched approximations. Future, ongoing work includes the enrichment of the smoothing space with suitable enrichment functions in order to investigate any additional accuracy improvements as well as the introduction of the approach in recently developed stable extended finite element schemes \cite{agathoschatzi2016,agathoschatzi2016a}.

\newpage

\bibliographystyle{plain}
\bibliography{main_Lsmxfem_v2}

\end{document}